\documentstyle[fleqn,12pt]{article}
\begin{document}\pagenumbering{arabic}\setcounter{page}{1}\pagestyle{plain}
\baselineskip=16pt

\thispagestyle{empty}
\rightline{MSUMB 97-09, December 1997} 
\vspace{1.4cm}

\begin{center}
{\Large\bf Hopf algebra structure of Gr$_q(1\vert 1)$ \\related to 
  GL$_q(1\vert 1)$} 
\end{center}

\vspace{1cm}
\noindent
Salih Celik \footnote{E-mail: scelik@fened.msu.edu.tr}
\footnote{{\bf E-mail}: sacelik@yildiz.edu.tr}\\ 
Mimar Sinan University, Department of Mathematics, \\
80690 Besiktas, Istanbul, TURKEY.

\vspace{1.5cm}
{\bf Abstract}

We show that the algebra of functions on the Grassmann supergroup 
Gr$_q(1\vert 1)$ has a (graded) Hopf algebra structure related to 
GL$_q(1\vert 1)$.

\vfill\eject
In the past few years, quantum groups$^1$ and $q$-deformed universal 
enveloping algebras$^2$ have been intensively studied both by 
mathematicians and mathematical physicists. From a mathematical point of 
view, these algebraic structures are just special classes of noncommutative 
Hopf algebras. 

The algebraic structure underlying quantum groups extends the theory of 
the supergroups.$^3$ The simplest quantum supergroup is GL$_q(1\vert 1)$, 
i.e. the deformation of the supergroup of 2x2 matrices with two 
bosonic (even) and two fermionic (odd) matrix entries. 

The aim of the present work is to construct the (graded) Hopf algebra 
structure of the Grassmann supergroup Gr$_q(1\vert 1)$, the superdual of 
GL$_q(1\vert 1)$, which was introduced in Ref. 4. Before discussing the 
(graded) Hopf algebra structure of Gr$_q(1\vert 1)$, let us first give 
some notations and useful formulas about the quantum Grassmann supergroup 
Gr$_q(1\vert 1)$. 

A Grassmann supermatrix $\hat{T}$ which is an element of Gr$(1\vert 1)$, 
is of the form 
$$ \hat{T} = \left(\matrix{ \alpha & b \cr c & \delta \cr} \right)$$
with two odd (greek letters) and two even (latin letters) matrix 
elements. The symbol {\it hat} is used to distinguish $\hat{T}$ from 
an element $T$ of GL$_q(1\vert 1)$. 

The $q$-deformation of the Grassman supergroup Gr$(1\vert 1)$ as a 
quantum matrix supergroup Gr$_q(1\vert 1)$ is generated by $\alpha$, 
$b$, $c$, $\delta$ with the relations$^4$ : 
$$ \alpha b = q^{-1} b \alpha, \qquad \alpha c = q^{-1} c \alpha, $$
$$ \delta b = q^{-1} b \delta, \qquad \delta c = q^{-1} c \delta, $$
$$\alpha \delta + \delta \alpha = 0, \qquad \alpha^2 = 0 = \delta^2, 
   \eqno(1)$$
$$ bc = cb + (q - q^{-1}) \delta \alpha $$
where $q$ is a non-zero complex number and $q^2 \neq 1$. The associative 
algebra (1) is equivalent to equation$^5$ 
$$R^1 \hat{T}_1 \hat{T}_2 = - \hat{T}_2 \hat{T}_1 R^2 \eqno(2)$$
where 
$$R^1 = \left(\matrix{ 
 q & 0          & 0   & 0 \cr 
 0 & - 1        & 0   & 0 \cr 
 0 & q - q^{-1} & - 1 & 0 \cr 
 0 & 0          & 0   & q^{-1} \cr }\right), \qquad 
R^2 = \left(\matrix{ 
 q^{-1} & 0   & 0          & 0 \cr 
 0      & - 1 & q^{-1} - q & 0 \cr 
 0      & 0   & - 1        & 0 \cr 
 0      & 0   & 0          & q \cr }\right), \eqno(3)$$
are both the solutions of the quantum (graded) Yang-Baxter equation. 
Here, we used the tensoring convention 
$$(\hat{T}_1)^{ij}_{kl} = (\hat{T} \otimes I)^{ij}_{kl} 
                        = \hat{T}^i_k \delta^j_l, $$
$$(\hat{T}_2)^{ij}_{kl} = (I \otimes \hat{T})^{ij}_{kl} 
                        = (-1)^{i(j+l)} \hat{T}^j_l \delta^i_k. \eqno(4)$$
The central element of the algebra (1) is$^4$ 
$$\hat{D}_q = b c^{-1} - \alpha c^{-1} \delta c^{-1} 
            = c^{-1} b - c^{-1} \alpha c^{-1} \delta. \eqno(5)$$
      
We now denote the algebra generated by the elements $\alpha$, $b$, $c$, 
$\delta$ with the relations (1) by $\hat{\cal A}$. We want to 
make the algebra $\hat{\cal A}$ into a (graded) Hopf algebra related to 
the quantum supergroup GL$_q(1\vert 1)$. Because of this, we state briefly 
some properties of the quantum supergroup GL$_q(1\vert 1)$ we are going to 
need in this work. 

The quantum supergroup GL$_q(1\vert 1)$ is generated by four generators 
$a$, $\beta$, $\gamma$, $d$ and the $q$-commutation relations $^3$ 
$$ a \beta = q \beta a, \qquad d \beta = q \beta d, $$
$$ a \gamma = q \gamma a, \qquad d \gamma = q \gamma d, \eqno(6)$$
$$ \beta \gamma + \gamma \beta = 0, \qquad \beta^2 = 0 = \gamma^2, $$
$$ a d = d a + (q - q^{-1}) \gamma \beta. $$
The generators satisfying the relations (6) generate the algebra 
called the algebra of functions on the quantum supergroup GL$_q(1\vert 1)$ 
and we shall denote it by ${\cal A}$. We know that the algebra 
${\cal A}$ is a (graded) Hopf algebra whose structure we now discuss. 
We represent the set of generators $a$, $\beta$, 
$\gamma$, $d$ in the form of a matrix 
$$T = \left(\matrix{ a & \beta \cr \gamma & d \cr} \right).$$ 
Then the relations (6) are equivalent to equation 
$$R T_1 T_2 = T_2 T_1 R \eqno(7)$$
where 
$$R = \left(\matrix{ 
 q & 0          & 0 & 0 \cr 
 0 & 1          & 0 & 0 \cr 
 0 & q - q^{-1} & 1 & 0 \cr 
 0 & 0          & 0 & q^{-1} \cr }\right). \eqno(8)$$
The superinverse of $T$ is given by $^3$ 
$$T^{-1} = \left(\matrix{ 
 A      & \Omega \cr 
 \Gamma & D \cr }\right) = 
 \left(\matrix{ 
 a^{-1} + a^{-1} \beta d^{-1} \gamma a^{-1} & - a^{-1} \beta d^{-1} \cr 
 - d^{-1} \gamma a^{-1} & d^{-1} + d^{-1} \gamma a^{-1} \beta d^{-1} \cr }
\right), \eqno(9)$$ 
and the superdeterminant is 
$$D_q = a d^{-1} - \beta d^{-1} \gamma. \eqno(10)$$
It is easy to verify that $D_q$ commutes with all matrix elements of $T$. 
Note that the matrix elements of $T$ with those of $T^{-1}$ satisfy the 
following relations 
$$a A = q^2 A a + 1 - q^2, \qquad d A = A d, $$
$$a D = D a, \qquad d D = q^2 D d + 1 - q^2, $$
$$a \Omega = q \Omega a, \qquad d \Omega = q \Omega d, $$
$$a \Gamma = q \Gamma a, \qquad d \Gamma = q \Gamma d, $$
$$\beta A = q A \beta, \qquad \gamma A = q A \gamma, \eqno(11)$$
$$\beta D = q D \beta, \qquad \gamma D = q D \gamma,$$
$$\beta \Omega = \Omega \beta, \qquad \gamma \Omega = - q^2 \Omega \gamma,$$
$$\beta \Gamma = - q^2 \Gamma \beta, \qquad \gamma \Gamma = \Gamma \gamma.$$

The usual coproduct is given by 
$$\Delta : {\cal A} \longrightarrow {\cal A} \otimes {\cal A}, \qquad 
   \Delta(t^i_j) = t^i_k \otimes t^k_j \eqno(12)$$
where summation over repeated indices is understood. One can rewrite the 
last formula in the following nice and elegant form 
$$\Delta(T) = T \dot{\otimes} T \eqno(13)$$
where $\otimes$ stands for the usual tensor product and the dot refers to the 
summation over repeated indices and reminds us about the usual matrix 
multiplication. 
The counit is given by 
$$ \varepsilon : {\cal A} \longrightarrow {\cal C}, \qquad 
    \varepsilon(t^i_j) = \delta^i_j. \eqno(14)$$
The coinverse (antipode) is given by 
$$S : {\cal A} \longrightarrow {\cal A}, \qquad S(T) = T^{-1}. \eqno(15)$$
It is not difficult to verify the following properties of the 
co-structures: 
$$(\Delta \otimes \mbox{id}) \circ \Delta = 
  (\mbox{id} \otimes \Delta) \circ \Delta, \eqno(16)$$
$$\mu \circ (\varepsilon \otimes \mbox{id}) \circ \Delta 
  = \mu' \circ (\mbox{id} \otimes \varepsilon) \circ \Delta, \eqno(17)$$
$$m \circ (S \otimes \mbox{id}) \circ \Delta = \varepsilon 
  = m \circ (\mbox{id} \otimes S) \circ \Delta, \eqno(18)$$
where id denotes the identity mapping, 
$$\mu : {\cal C} \otimes {\cal A} \longrightarrow {\cal A}, \qquad 
  \mu' : {\cal A} \otimes {\cal C} \longrightarrow {\cal A}$$
are the canonical isomorphisms, defined by 
$$\mu(k \otimes a) = ka = \mu'(a \otimes k), \quad \forall a \in {\cal A}, 
  \quad \forall k \in {\cal C} \eqno(20)$$
and $m$ is the multiplication map 
$$m : {\cal A} \otimes {\cal A} \longrightarrow {\cal A}, \qquad 
  m(a \otimes b) = ab. \eqno(21)$$
The multiplication in ${\cal A} \otimes {\cal A}$ follows the rule 
$$(A \otimes B) (C \otimes D) = (-1)^{p(B) p(C)} AC \otimes BC \eqno(22)$$
where $p(X)$ is the $z_2$-grade of $X$, i.e. $p(X) = 0$ for even variables 
and $p(X) = 1$ for odd variables. 

Since the (graded) Hopf algebra structure of Gr$_q(1\vert 1)$ is related 
to those of GL$_q(1\vert 1)$ it is necessary to obtain the commutation 
relations of the generators of $\hat{\cal A}$ with those of ${\cal A}$. 
We define the (mixed) commutation relations between the generators of 
$\hat{\cal A}$ and ${\cal A}$ as follows: 
$$R \hat{T}_1 T_2 = (-1)^{p(T_2)} T_2 \hat{T}_1 R', \qquad 
  R' = R - (q - q^{-1}) P \eqno(23)$$
where $P$ is the superpermutation matrix. The equation (23) gives the mixed 
relations 
$$a \alpha = q^2 \alpha a, \qquad \beta \alpha = - q \alpha \beta, $$
$$a b = q b a + (q^2 - 1) \alpha \beta, \qquad \beta b = b \beta, $$
$$a c = q c a + (q^2 - 1) \alpha \gamma, \qquad 
  \beta c = c \beta + (q - q^{-1}) \alpha d,$$
$$a \delta = \delta a + (q - q^{-1}) (\beta c - b \gamma), \qquad 
  \beta \delta = - q^{-1} \delta \beta + (1 - q^{-2}) b d, \eqno(24)$$
$$d \alpha = \alpha d, \qquad \gamma \alpha = - q \alpha \gamma, $$
$$d b = q^{-1} b d, \qquad \gamma b = b \gamma - (q - q^{-1}) \alpha d, $$
$$d c = q^{-1} c d, \qquad \gamma c = c \gamma,$$
$$d \delta = q^{-2} \delta d, \qquad 
  \gamma \delta = - q^{-1} \delta \gamma + (1 - q^{-2}) c d. $$
Using these relations, it is easy to verify that $\hat{D}_q$, which is given 
by (5), is still a central element, i.e. $\hat{D}_q$ also commutes with 
the generators of ${\cal A}$. 

After some algebra the commutation relations between the matrix elements of 
$\hat{T}$ with $T^{-1}$ are obtained to be  
$$\alpha A = q^2 A \alpha, \qquad \delta A = A \delta,$$
$$\alpha D = D \alpha, \qquad \delta D = q^{-2} D \delta + 
  (q - q^{-1})^2 A \alpha + (q^{-2} - 1) (\Omega c - \Gamma b),$$
$$\alpha \Omega = - q \Omega \alpha, \qquad 
  \delta \Omega = - q^{-1} \Omega \delta + (q^{-1} - q) A b,$$
$$\alpha \Gamma = - q\Gamma \alpha, \qquad 
  \delta \Gamma = - q^{-1} \Gamma \delta + (q^{-1} - q) A c, \eqno(25)$$
$$b A = q A b, \qquad c A = q A c,$$
$$b D = q^{-1} D b + (q - q^{-1}) \Omega \alpha, \qquad 
  c D = q^{-1} D c + (q - q^{-1}) \Gamma \alpha,$$
$$b \Omega = \Omega b, \qquad c \Omega = \Omega c + (q^2 - 1) A \alpha,$$
$$b \Gamma = \Gamma b + (1 - q^2) A \alpha, \qquad 
  c \Gamma = \Gamma c,$$
and 
$$D_q u = q^2 u D_q, \qquad u \in \{\alpha, b, c, \delta\}.\eqno(26)$$

Before defining a coproduct on the algebra $\hat{\cal A}$, let us note the 
following facts. Let $\hat T$ and $\hat{T}'$ be any two supercommuting 
matrices that satisfy (1). We denote a product $\hat{T} \hat{T}'$ by $T$. 
Then, it can be verified that the matrix elements of $T$ satisfy the 
commutation relations (6) of GL$_q(1\vert 1)$, i.e. if 
$$T = \left(\matrix{ \alpha & b \cr c & \delta \cr}\right) 
      \left(\matrix{ \alpha' & b' \cr c' & \delta' \cr}\right) $$
then we have the relations (6). In short, 
$$\hat{T}, \hat{T}' \in Gr_q(1\vert 1) ~\Longrightarrow~ 
  T = \hat{T} \hat{T}' \in GL_q(1\vert 1). $$
In view of these facts, we can say that there may be no coproduct of 
the usual form $\Delta(\hat{T}) = \hat{T} \dot{\otimes} \hat{T}$. For, this 
coproduct, if existed, would be invariant under the $q$-commutation 
relations (6) of GL$_q(1\vert 1)$. But we can define a map on the 
algebra $\hat{\cal A}$ as follows: 
$$\hat{\Delta} : \hat{\cal A} \longrightarrow \hat{\cal A} \otimes 
  \hat{\cal A}, \qquad 
   \hat{\Delta}(\hat{T}) = \hat{T} \dot{\otimes} T + 
   (-1)^{p(T)} T \dot{\otimes} \hat{T}. 
   \eqno(27)$$
Explicity, 
$$\hat{\Delta}(\alpha) = \alpha \otimes a + b \otimes \gamma + 
   a \otimes \alpha - \beta \otimes c, $$ 
$$\hat{\Delta}(b) = b \otimes d + \alpha \otimes \beta + 
   a \otimes b - \beta \otimes \delta, \eqno(28)$$ 
$$\hat{\Delta}(c) = c \otimes a + \delta \otimes \gamma - 
   \gamma \otimes \alpha + d \otimes c, $$ 
$$\hat{\Delta}(\delta) = \delta \otimes d + c \otimes \beta - 
   \gamma \otimes b + d \otimes \delta. $$ 

The action on the generators of $\hat{\cal A}$ of 
$\hat{\varepsilon} : \hat{\cal A} \longrightarrow {\cal C}$ is 
$$\hat{\varepsilon}(\alpha) = \hat{\varepsilon}(b) = 
   \hat{\varepsilon}(c) = \hat{\varepsilon}(\delta) = 0. \eqno(29)$$
Finally, we define the coinverse as 
$$\hat{S} : \hat{\cal A} \longrightarrow \hat{\cal A}, \qquad 
  \hat{S}(\hat{T}) = - (-1)^{p(T^{-1})} T^{-1} \hat{T} T^{-1}. \eqno(30)$$
The action of $\hat{S}$ on the generators of $\hat{\cal A}$ is 
$$\hat{S}(\alpha) = - (\alpha A + b \Gamma) A + 
   q (c A + \delta \Gamma) \Omega, $$
$$\hat{S}(b) = - (\alpha A - c \Omega) \Omega - 
  q (b A + \delta \Omega) D, \eqno(31)$$
$$\hat{S}(c) = - (\alpha A + b \Gamma) \Gamma - q (c A + \delta \Gamma) D, $$
$$\hat{S}(\delta) = q^2 (\alpha A - c \Omega) D + 
  q (\alpha \Omega + q^2 b D) \Gamma - (\alpha A + q^2 \delta D) D. $$
It is not difficult to check that the maps $\hat{\Delta}$ and 
$\hat{\varepsilon}$ are both algebra homomorphisms and $\hat{S}$ is an 
algebra anti-homomorphism and also the three maps satisfy the properties 
(16)-(18), and they preserve the relations (24) provided that the 
action on the generators of ${\cal A}$ of $\hat{\Delta}$ is the same 
with (13). 

The coproduct, counit and coinverse which are specified above supply 
Gr$_q(1\vert 1)$ with a structure, which can be called a quasi-Hopf algebra. 

It is interesting to note that there is a close connection with the 
differential calculus$^6$ on the quantum supergroup GL$_q(1\vert 1)$ 
via the equation (23). In fact we have observed that the matrix elements 
of $\hat T \in \mbox{Gr}_q(1\vert 1)$ are just the differentials of the 
matrix elements of $T \in \mbox{GL}_q(1\vert 1)$. In other words, we can 
interpret the generating elements of Gr$_q(1\vert 1)$ as differentials 
of coordinate functions on GL$_q(1\vert 1)$. In this case, we can write 
$\hat{T} = {\sf d} T$ (more information on these issues are given in ref. 6). 
Then the extended algebra can be interpreted as an algebra of differential 
forms on GL$_q(1\vert 1)$. Thus the coproduct is interpreted as a (left and 
right) coaction of the quantum supergroup GL$_q(1\vert 1)$ on differential 
forms. To this end, we consider the two maps 
$$\Delta_R : \Gamma \longrightarrow \Gamma \otimes {\cal A}, \qquad 
 \Delta_R \circ {\sf d} = 
 ({\sf d} \otimes \mbox{id}) \circ \Delta \eqno(32\mbox{a})$$
and 
$$\Delta_L : \Gamma \longrightarrow {\cal A} \otimes \Gamma, \qquad 
  \Delta_L \circ {\sf d} = (\tau \otimes {\sf d}) \circ \Delta, 
  \eqno(32\mbox{b})$$
where $\Gamma$ denotes the differential algebra of ${\cal A}$. 
Here $\tau: \Gamma \longrightarrow \Gamma$ is the linear map of degree zero 
which gives $\tau(a) = (-1)^{p(a)} a$. 
We now define a map $\phi_R$ as follows 
$$\phi_R(u_1 {\sf d}v_1 + {\sf d}v_2 u_2) = 
 \Delta(u_1) \Delta_R({\sf d}v_1) + \Delta_R({\sf d}v_2) \Delta(u_2)\eqno(33)$$
and another map $\phi_L$ by replacing $L$ with $R$. 
The following identities are satisfied 
$$(\phi_R \otimes \mbox{id}) \circ \phi_R = (\mbox{id} \otimes \Delta) \circ \phi_R 
  \qquad (\mbox{id} \otimes \epsilon) \circ \phi_R = \mbox{id}, \eqno(34\mbox{a})$$
and 
$$(\mbox{id} \otimes \phi_L) \circ \phi_L = (\Delta \otimes \mbox{id}) \circ \phi_L 
  \qquad (\epsilon \otimes \mbox{id}) \circ \phi_L = \mbox{id}. \eqno(34\mbox{b})$$
Consequently, we define the map $\hat{\Delta}$, in (27), as 
$$\hat{\Delta} = \phi_R + \phi_L. \eqno(35)$$

\noindent
{\bf ACKNOWLEDGEMENT}

This work was supported in part by {\bf T. B. T. A. K.} the 
Turkish Scientific and Technical Research Council. 

I would like to express my deep gratitude to the referee for critical 
comments and suggestions on the manuscript. 

\vfill\eject
\noindent
$^1$ N. Y. Reshetikhin, L. A. Takhtajan and L. D. Faddeev, 
    Leningrad Math. J. {\bf 1} \\
\hspace*{0.2cm} (1990), 193. \\
$^2$ V. G. Drinfeld, Quantum groups, 
    in {\it Proc. } IMS, Berkeley, (1986). \\
$^3$ Yu. I. Manin, Commun. Math. Phys. {\bf 123} (1989), 163; \\
\hspace*{0.2cm} E. Corrigan, D. Fairlie, P. Fletcher and R. Sasaki, 
 J. Math. Phys. {\bf 31} (1990), 
\hspace*{0.2cm} 776;\\
\hspace*{0.2cm} S. Schwenk, B. Schmidke and S. Vokos, Z. Phys. C {\bf 46} (1990), 643; \\
\hspace*{0.2cm} B. Schmidke, S. Vokos, and B. Zumino, Z. Phys. C {\bf 48} (1990), 249. \\
$^4$ S. \c Celik and S. A. \c Celik, Balkan Phys. Lett. {\bf 3} (1995), 188. \\
$^5$ S. \c Celik, J. Math. Phys. {\bf 37} (1996), 3568; \\
\hspace*{0.2cm} S. \c Celik, Balkan Phys. Lett. {\bf 5} (1997), 149. \\
$^6$ S. \c Celik and S. A. \c Celik, "On the differential geometry of 
 GL$_q(1\vert 1)$", J. Phys. \\
\hspace*{0.2cm} A: Math. Gen. (in press) (1998).

\end{document}